\documentclass[10pt]{amsart}
\usepackage{amsmath, amsfonts, amssymb, amstext, amsthm, graphicx, subfigure}
\newtheorem{rotation}{Proposition}
\newtheorem{rotation2}[rotation]{Corollary}
\newtheorem{basis}[rotation]{Lemma}
\newtheorem{decomposition}[rotation]{Lemma}
\newtheorem{distance}[rotation]{Lemma}

\begin{document}

\title[Area of Revolution about a Line]{The Area of the Surface Generated by Revolving a Graph About Any Line}

\author{Edray Herber Goins}
\address{Purdue University \\ Department of Mathematics \\ 150 North University Street \\ West Lafayette, IN 47907}
\email{egoins@math.purdue.edu}

\author{Talitha M. Washington}
%\address{Department of Mathematics \\ 1800 Lincoln Avenue \\ University of Evansville \\ Evansville, IN 47722}
\address{Howard University \\ College of Arts and Sciences \\ Department of Mathematics \\ 204 Academic Support Building B \\ Washington, DC 20059}
\email{tw65@evansville.edu}



\thanks{This exposition is dedicated to the inquisitive students in advanced Calculus courses who enjoy seeing mathematics beyond those topics covered in a standard textbook on the subject.}

\begin{abstract}
We discuss a general formula for the area of the surface that is generated by a graph $[t_0, t_1] \to \mathbb R^2$ sending $t \mapsto \bigl( x(t), \, y(t) \bigr)$ revolved around a general line $L: \, A \, x + B \, y = C$.  As a corollary, we obtain a formula for the area of the surface formed by revolving $y = f(x)$ around the line $y = m \, x + k$. 
\end{abstract}

\maketitle

\section{Introduction}

\begin{quote} 
\textit{You spin me right round, baby, right round! \\ -- Dead or Alive}
\end{quote}
\vskip 0.2in

Many students in advanced Calculus courses learn how to compute the area of a surface of revolution.  Perhaps the best known example goes something like this: For functions $y = f(x)$ differentiable on the interval $a \leq x \leq b$, the area of the surface of revolution of its graph about the $x$-axis is given by the integral
\begin{equation} \label{x-axis} 2 \, \pi \int_a^b |y| \, \sqrt{1 + \left( \dfrac {dy}{dx} \right)^2} \, dx = 2 \, \pi \int_a^b \bigl| f(x) \bigr| \, \sqrt{1 + \bigl[ f'(x) \bigr]^2} \, dx. \end{equation}

\noindent An example of this type of rotation can be found in Figure \ref{rotation-x} below.  

\begin{center} \begin{figure}[h] 
\caption{Surfaces of Revolution of $y = x^2 - 3 \, x + 12$} \label{rotation-x}
\subfigure[Graph of Function]{\includegraphics[height=0.23\textheight]{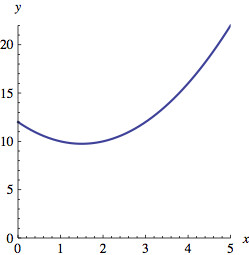}} \quad
\subfigure[Revolved around $x$-axis]{\includegraphics[height=0.23\textheight]{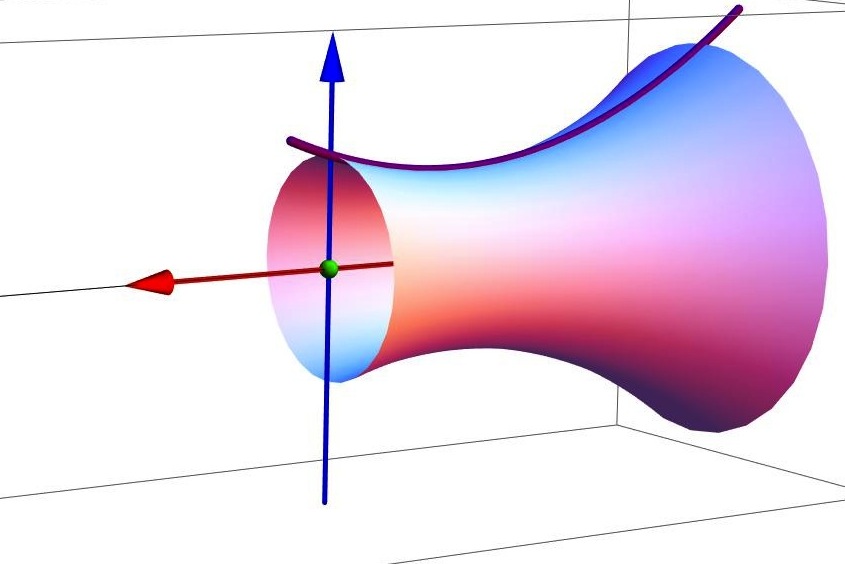}} \quad
\end{figure} \end{center}

\noindent There is a similar formula for computing the area of the surface of revolution about the $y$-axis, but one must work with the inverse function $x = g(y)$ and assume that it is a function differentiable on the interval $c \leq y \leq d$:
\begin{equation} \label{y-axis} 2 \, \pi \int_c^d |x| \, \sqrt{1 + \left( \dfrac {dx}{dy} \right)^2} \, dy = 2 \, \pi \int_c^d \bigl| g(y) \bigr| \, \sqrt{1+ \bigl[ g'(y) \bigr]^2} \, dy. \end{equation}

\noindent As above, an example of this type of rotation can be found in Figure \ref{rotation-y} below.  

\begin{center} \begin{figure}[h] 
\caption{Surfaces of Revolution of $x = y^3 - 12 \, y^2 + 444 \, y + 62$} \label{rotation-y}
\subfigure[Graph of Function]{\includegraphics[height=0.30\textheight]{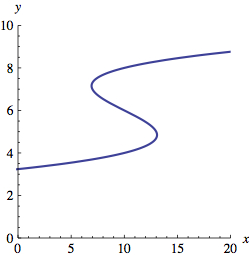}} \quad
\subfigure[Revolved around $y$-axis]{\includegraphics[height=0.30\textheight]{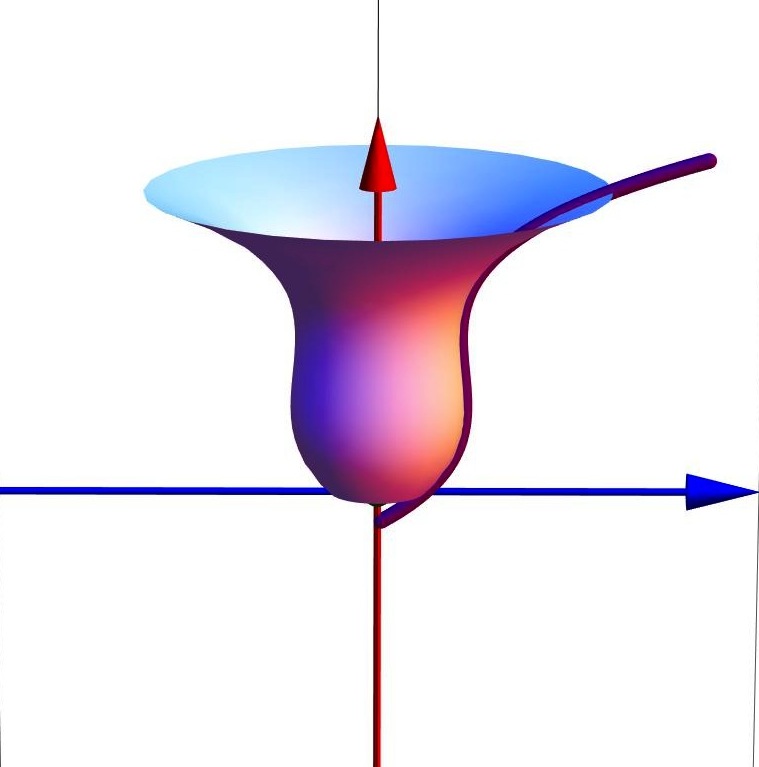}}
\end{figure} \end{center}

These are special cases of graphs which can be parametrized by continuous functions, say $x = x(t)$ and $y = y(t)$, which are both differentiable on the interval $t_0 \leq t \leq t_1$.  The area of the surface of revolution of parametrized graph is given by the integral
\begin{equation} \begin{aligned}
2 \, \pi \int_{t_0}^{t_1} |y(t)| \, \sqrt{ \left( \frac {dx}{dt} \right)^2 + \left( \frac {dy}{dt} \right)^2} \, dt &  & \text{when revolved about the $x$-axis,} \\[5pt]
2 \, \pi \int_{t_0}^{t_1} |x(t)| \, \sqrt{ \left( \frac {dx}{dt} \right)^2 + \left( \frac {dy}{dt} \right)^2} \, dt & & \text{when revolved about the $y$-axis.}
\end{aligned}\end{equation}

\noindent Such formulas can be found in Stewart's \textit{Calculus} \cite[Section 8.2, Section 10.2]{Stewart2008Calculus:-Early}.  A general discussion of integrals over surfaces can be found in Apostol's \textit{Calculus} \cite[Chapter 12]{MR0248290}.

\section{Rotations on a Slant}

These seemingly unrelated formulas must be intimidating for advanced Calculus students.  Which one should they use to solve problems?  Why are the formulas so different?  And where do these formulas come from in the first place?  As educators, it's our job to ease fears and provide the simplest explanations possible to even the most complex mathematical concepts.  In this note, we'll provide a relatively simple proof of the following generalization.

\begin{rotation} \label{rotation} Consider the graph parametrized by functions $x = x(t)$ and $y = y(t)$, both differentiable on the interval $t_0 \leq t \leq t_1$.  The area of the surface of revolution of this graph, when revolved about a line $L: \, A \, x + B \, y = C$, is given by the integral
\begin{equation} \begin{aligned} 2 \, \pi & \int_{t_0}^{t_1} r(t) \, \sqrt{ \left( \frac {dx}{dt} \right)^2 + \left( \frac {dy}{dt} \right)^2} \, dt \\ & = 2 \, \pi \int_{t_0}^{t_1} r(t) \, \sqrt{ \bigl[ x'(t) \bigr]^2 + \bigl[ y'(t) \bigr]^2} \, dt \end{aligned} \end{equation}
\noindent in terms of the nonnegative function
\begin{equation} r(t) = \frac {\bigl| A \, x(t) + B \, y(t) - C \bigr|}{\sqrt{A^2 + B^2}}. \end{equation} \end{rotation}

Before we present the proof, we explain how to use this generalization to derive Equations \eqref{x-axis} and \eqref{y-axis}.  For the former, the graph of $y = f(x)$ can be parametrized by the functions $x(t) = t$ and $y(t) = f(t)$ on the interval $a \leq t \leq b$.  The $x$-axis is simply the line $L: \, y = 0$, so we may choose $(A,B,C) = (0,1,0)$.  Then $r(t) = |f(t)|$, and we find Equation \eqref{x-axis}.  For the latter, the graph of $x = g(y)$ can be parametrized by the functions $x(t) = g(t)$ and $y(t) = t$ on the interval $c \leq t \leq d$.  The $y$-axis is simply the line $L: \, x = 0$, so we may choose $(A,B,C) = (1,0,0)$.  Then $r(t) = |g(t)|$, and we find Equation \eqref{y-axis}.

On pages 551 and 552 of Stewart's \textit{Calculus} \cite{Stewart2008Calculus:-Early}, there is a Discovery Project entitled ``Rotations on a Slant'', where a differentiable function $y = f(x)$ is rotated about the line $y = m \, x + k$ expressed in point-slope form.  As a solution to Exercise 5, we recover the following result.

\begin{rotation2} Consider a function $y = f(x)$ which is differentiable on the interval $a \leq x \leq b$.  The area of the surface of revolution of its graph about the line $y = m \, x + k$ is given by the integral
\begin{equation} \begin{aligned} 2 \, \pi & \int_a^b |y- m \, x - k| \, \sqrt{ \dfrac {1 + (dy/dx)^2}{1 + m^2}} \, dx \\ & = 2 \, \pi \int_a^b \bigl| f(x) - m \, x - k \bigr| \, \sqrt{\dfrac {1 + \bigl[ f'(x) \bigr]^2}{1 + m^2}} \, dx. \end{aligned} \end{equation} \end{rotation2}

We explain why this corollary is true.  As above, write $y = f(x)$ in terms of the functions $x(t) = t$ and $y(t) = f(t)$ on the interval $a \leq t \leq b$.  The line $L: \, y = m \, x + k$ can be expressed in the form $A \, x + B \, y = C$ upon choosing $(A,B,C) = (-m,1,k)$, so that we have the nonnegative function
\begin{equation} r(t) = \frac {\bigl| A \, x(t) + B \, y(t) - C \bigr|}{\sqrt{A^2 + B^2}} = \dfrac {\bigl| f(t) - m \, t - k \bigr|}{\sqrt{1 + m^2}}. \end{equation}

\noindent Proposition \ref{rotation} implies that we have the area
\begin{equation} \begin{aligned} 2 \, \pi & \int_{t_0}^{t_1} r(t) \, \sqrt{ \bigl[ x'(t) \bigr]^2 + \bigl[ y'(t) \bigr]^2} \, dt \\ & = 2 \, \pi \int_a^b \dfrac {\bigl| f(x) - m \, x - k \bigr|}{\sqrt{1 + m^2}} \sqrt{1 + \bigl[ f'(x) \bigr]^2} \, dx. \end{aligned} \end{equation}

\noindent Some examples of this result can be found in Figure \ref{proposition} below.

\begin{center} \begin{figure}[h] 
\caption{Revolutions of $y = x^2 - 3 \, x + 12$ Around Lines} 
\label{proposition}
\subfigure[$y = 0$]{\includegraphics[height=0.19\textheight]{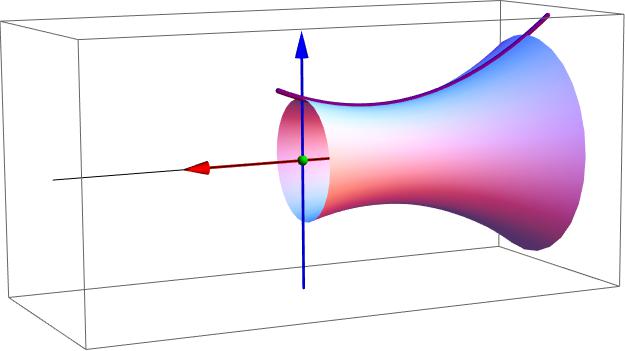}} 
\subfigure[$3 \, x + 4 \, y = 0$]{\includegraphics[height=0.19\textheight]{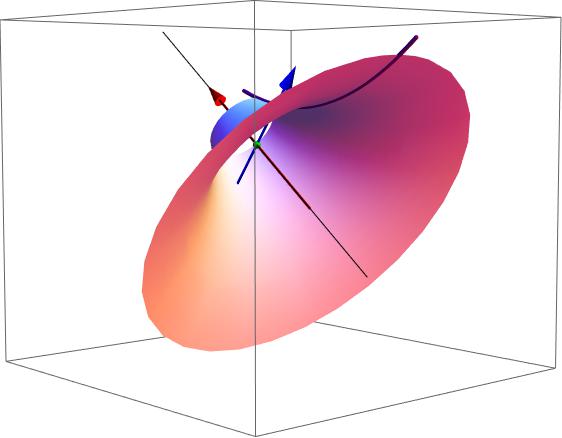}} 
\subfigure[$x = 0$]{\includegraphics[height=0.25\textheight]{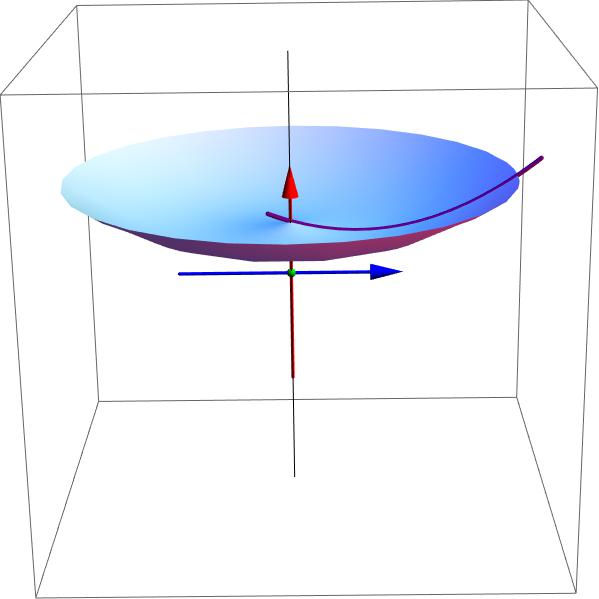}} 
\subfigure[$3 \, x - 4 \, y = 0$]{\includegraphics[height=0.25\textheight]{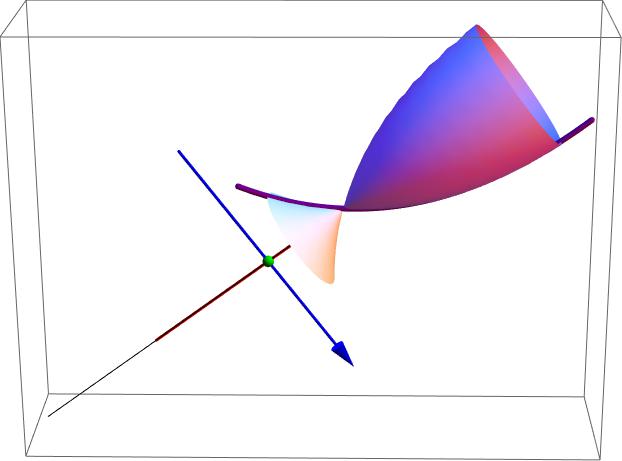}}
\end{figure} \end{center}

\section{Proof of Main Result}

The idea for the proof will be to create a new coordinate system using the line $L$ as an axis.   To this end, let's define the following three points in the plane:
\begin{equation} \begin{aligned} \mathbf{O} & = \left( \frac {A \, C}{A^2 + B^2}, \ \frac {B \, C}{A^2 + B^2} \right), \\[5pt]
\mathbf{u} & = \left( - \frac {B}{\sqrt{A^2 + B^2}}, \ \frac {A}{\sqrt{A^2 + B^2}} \right), \\[5pt]
\mathbf{v} & = \left( \frac {A}{\sqrt{A^2 + B^2}}, \ \frac {B}{\sqrt{A^2 + B^2}} \right).  
\end{aligned} \end{equation}

\begin{center} \begin{figure}[h] 
\caption{Plot of $\mathbf{0}$, $\mathbf{u}$, and $\mathbf{v}$ on the Line $L: \, A \, x + B \, y = C$} \label{axes}
\includegraphics[width=0.8\textwidth]{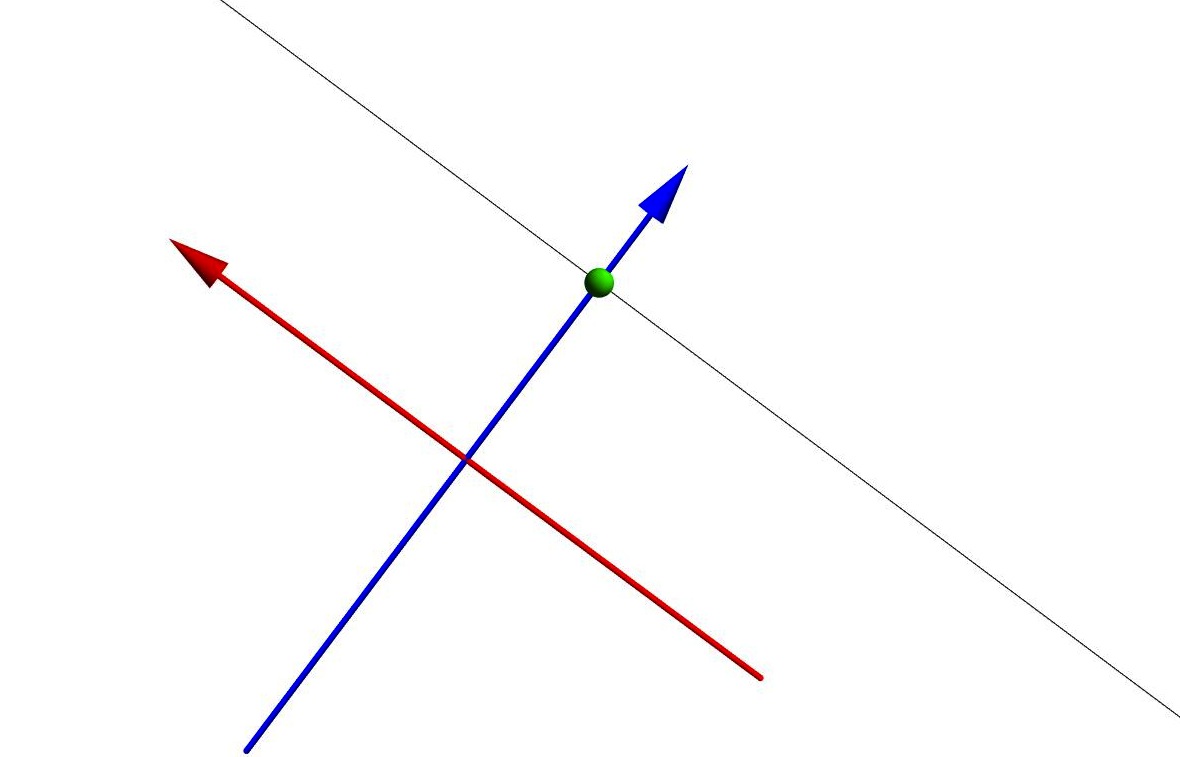} 
\end{figure} \end{center}

\noindent Each of these points has a geometric interpretation.  The first point $\mathbf{O}$ will act as our ``origin''; it lies on the line $L$.  The second point $\mathbf{u}$ lies along the direction of $L$, while the third point $\mathbf{v}$ lies perpendicular to $L$.  A diagram of these points with respect to $L$ can be found in Figure \ref{axes}.

In order to prove Proposition \ref{rotation}, we will prove a series of results.

\begin{basis} \label{basis} $\{ \mathbf{u}, \, \mathbf{v} \}$ forms an orthonormal basis for $\mathbb R^2$. \end{basis}

Recall that an orthonormal basis $\{ \mathbf{u}, \, \mathbf{v} \}$ is a set of two vectors which are orthogonal (that is, $\mathbf{u} \cdot \mathbf{v} = 0$) and normal (that is, $\| \mathbf{u} \| = \| \mathbf{v} \| = 1$).  This is indeed the case for the vectors above because we have the identities
\begin{equation} \begin{aligned}
\mathbf{u} \cdot \mathbf{v} & = - \frac {B}{\sqrt{A^2 + B^2}} \ \frac {A}{\sqrt{A^2 + B^2}} + \frac {A}{\sqrt{A^2 + B^2}} \ \frac {B}{\sqrt{A^2 + B^2}} \\ & = \dfrac {-B \, A + A \, B}{A^2 + B^2} = 0, \\[5pt]
\| \mathbf{u} \| & = \sqrt{ \left( - \frac {B}{\sqrt{A^2 + B^2}} \right)^2  + \left( \frac {A}{\sqrt{A^2 + B^2}} \right)^2} = \sqrt{\dfrac {B^2 + A^2}{A^2 + B^2}} = 1, \\[5pt]
\| \mathbf{v} \| & = \sqrt{ \left( \frac {A}{\sqrt{A^2 + B^2}} \right)^2 + \left( - \frac {B}{\sqrt{A^2 + B^2}} \right)^2} = \sqrt{\dfrac {A^2 + B^2}{A^2+ B^2}} = 1.
\end{aligned} \end{equation}

Next, we'll use these three points $\mathbf{O}$, $\mathbf{u}$, and $\mathbf{v}$ to decompose points $(x,y)$ in the plane along the axes defined above.

\begin{decomposition} \label{decomposition} Every point in the plane can be decomposed into parts that lie on the line $L$ and parts that lie perpendicular to the line.  Explicitly,
\begin{equation} \bigl( x, y \bigr) = \underbrace{\left[ \mathbf{O} + \frac {-B \, x + A \, y}{\sqrt{A^2 + B^2}} \ \mathbf{u} \right]}_{\text{on the line $L$}} + \underbrace{\left[ \frac {A \, x + B \, y - C}{\sqrt{A^2 + B^2}} \ \mathbf{v} \right]}_{\text{perpendicular to the line $L$}}. \end{equation}
\end{decomposition}

In order to prove this, we must perform some algebra.  The decomposition holds because we have the sum
\begin{equation} \begin{tabular}{rccccc} 
$\mathbf{O}$ & = & $\biggl($ & $\dfrac {\quad \qquad \qquad \qquad A \, C}{A^2 + B^2}$, & $\dfrac {\quad \qquad \qquad \qquad B \, C}{A^2 + B^2}$ & $\biggr)$ \\[10pt]
$\dfrac {-B \, x + A \, y}{\sqrt{A^2 + B^2}} \ \mathbf{u}$ & = & $\biggl($ & $\dfrac {B^2 \, x - A \, B \, y \quad \qquad}{A^2 + B^2}$, & $\dfrac {- A \, B \, x + A^2 \, y \quad \qquad}{A^2 + B^2}$ & $\biggr)$ \\[10pt]
$\dfrac {A \, x + B \, y - C}{\sqrt{A^2 + B^2}} \ \mathbf{v}$ & = & $\biggl($ & $\dfrac {A^2 \, x + A \, B \, y - A \, C}{A^2 + B^2}$, & $\dfrac {A \, B \, x + B^2 \, y - B \, C}{A^2 + B^2}$ & $\biggr)$ \\[10pt] \hline
& & & & & \\[-8pt]
\text{Sum of Points} & = & $\biggl($ & $\dfrac {A^2 + B^2}{A^2 + B^2} \, x$, & $\dfrac {A^2 + B^2}{A^2 + B^2} \, y$ & $\biggr)$
\end{tabular} \end{equation}

\noindent Moreover, the sum of the first two points
\begin{equation}
\mathbf{O} + \frac {-B \, x + A \, y}{\sqrt{A^2 + B^2}} \ \mathbf{u} = \left( \dfrac{ B^2 \, x - A \, B \, y + A \, C}{A^2 + B^2}, \ \dfrac {-A \, B \, x + A^2 \, y + B \, C}{A^2 + B^2} \right)
\end{equation}

\noindent lies on the line $L$:
\begin{equation} \begin{aligned}
A & \left( \dfrac{ B^2 \, x - A \, B \, y + A \, C}{A^2 + B^2} \right) + B \left( \dfrac {-A \, B \, x + A^2 \, y + B \, C}{A^2 + B^2} \right) \\ 
& = \dfrac{ \bigl( A \, B^2 \, x - A^2 \, B \, y + A^2 \, C \bigr)  + \bigl( - A \, B^2 \, x + A^2 \, B \, y + B^2 \, C \bigr)}{A^2 + B^2} = C.
\end{aligned} \end{equation}

\noindent The third point is perpendicular to the line because $\mathbf{v}$ is orthogonal to $\mathbf{u}$.

Next, we put these results together to compute a useful formula.

\begin{distance} \label{distance} The distance from a point $(x,y)$ to the line $L$ is given by
\begin{equation} r = \frac {\bigl| A \, x + B \, y - C \bigr|}{\sqrt{A^2 + B^2}}. \end{equation}
\end{distance}

To prove this, we essentially use the Pythagorean Theorem: we project any point $(x,y)$ in the plane into two directions, one along the line $L$ and the other perpendicular to it.  The distance from this point to the line is simply the length of the perpendicular direction, that is, the length of the vector along $\mathbf{v}$.
\begin{equation} r = \left \| \frac {A \, x + B \, y - C}{\sqrt{A^2 + B^2}} \ \mathbf{v} \right \| = \frac {\bigl| A \, x + B \, y - C \bigr|}{\sqrt{A^2 + B^2}} \, \bigl \| \mathbf{v} \bigr \| = \frac {\bigl| A \, x + B \, y - C \bigr|}{\sqrt{A^2 + B^2}}. \end{equation}

\textit{Proof of Proposition \ref{rotation}:}  Say that we have two functions $x = x(t)$ and $y = y(t)$, both differentiable on the interval $t_0 \leq t \leq t_1$, as well as a line $L: \, A \, x + B \, y = C$.  According to Lemma \ref{distance}, the distance from a point $(x,y) = \bigl( x(t), \, y(t) \bigr)$ to the line $L$ is given by the function 
\begin{equation} r(t) = \frac {\bigl| A \, x(t) + B \, y(t) - C \bigr|}{\sqrt{A^2 + B^2}}. \end{equation}

\noindent This will act as the radius of rotation around the line.  A small arc length along the graph at this point is given by the differential
\begin{equation}
ds = \sqrt{{dx}^2 + {dy}^2} = \sqrt{ \left( \frac {dx}{dt} \right)^2 + \left( \frac {dy}{dt} \right)^2} \, dt = \sqrt{ \bigl[ x'(t) \bigr]^2 + \bigl[ y'(t) \bigr]^2} \, dt. 
\end{equation}

\begin{center} \begin{figure}[h] \caption{Cylindrical Shell of Radius $r(t)$ and Length $ds$} \label{shell} \includegraphics[width=0.7\textwidth]{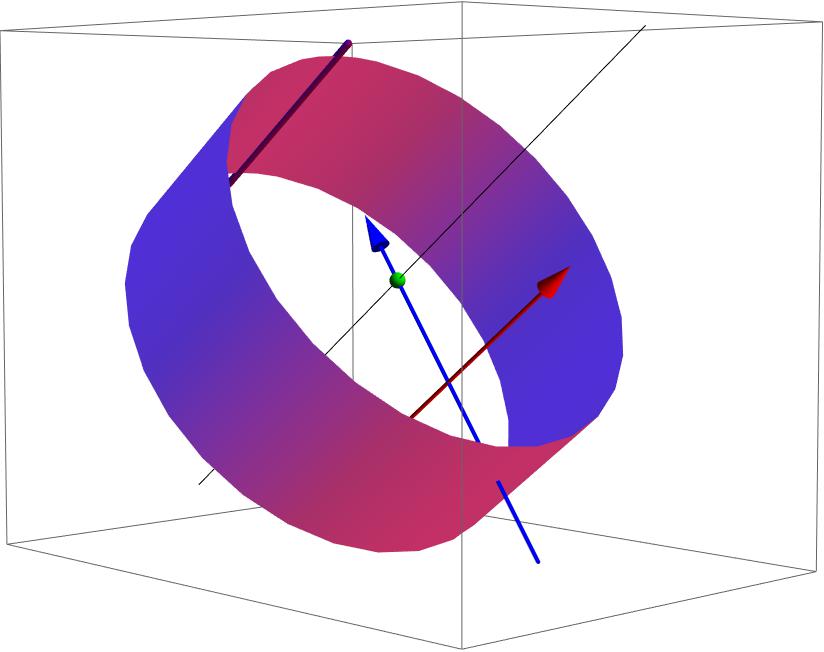} \end{figure} \end{center}

\noindent This will act as the width of the cylindrical shell.  (See Figure \ref{shell} for an illustration.) Putting these together, the area differential of the surface of revolution is  
\begin{equation} \begin{aligned} d \text{Area} & = 2 \, \pi \, r(t) \, ds \\ & = 2 \, \pi \, \frac {\bigl| A \, x(t) + B \, y(t) - C \bigr|}{\sqrt{A^2 + B^2}} \, \sqrt{ \bigl[ x'(t) \bigr]^2 + \bigl[ y'(t) \bigr]^2} \, dt, \end{aligned} \end{equation}

\noindent so that the desired area of revolution is given by
\begin{equation} \begin{aligned} \text{Area} & = 2 \, \pi \int_{t_0}^{t_1} r(t) \, d s \\ & = 2 \, \pi \int_{t_0}^{t_1} \frac {\bigl| A \, x(t) + B \, y(t) - C \bigr|}{\sqrt{A^2 + B^2}} \, \sqrt{ \bigl[ x'(t) \bigr]^2 + \bigl[ y'(t) \bigr]^2} \, dt. \end{aligned} \end{equation}

\section{Example}

\begin{center} \begin{figure}[h] \caption{Circle $x^2 + y^2 = r^2$ Revolved Around $A \, x + B \, y = C$} \label{torus} \includegraphics[width=0.7\textwidth]{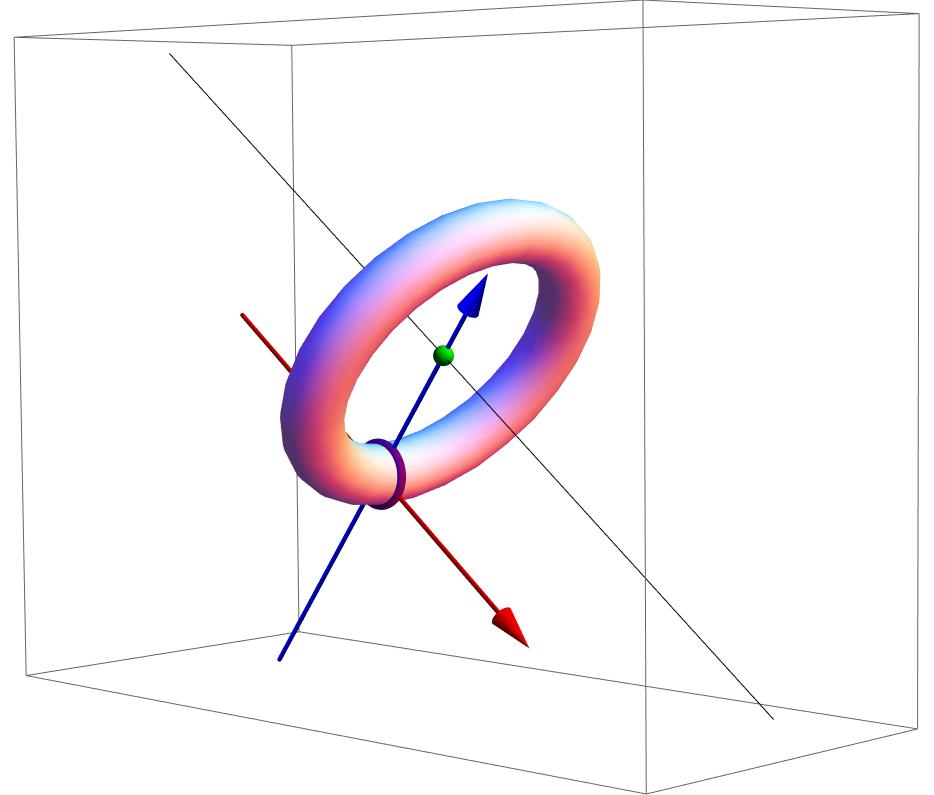} \end{figure} \end{center}

We illustrate these ideas with an example.  Say that we wish to rotate the circle $x^2 + y^2 = r^2$ of radius $r$ around the line $A \, x + B \, y = C$ -- but we'll assume for simplicity that
\begin{equation} R = \dfrac {C}{\sqrt{A^2 + B^2}} >  r.  \end{equation}
\noindent As shown in Figure \ref{torus}, this is just a torus.  We'll show that the area of this surface is $4 \, \pi^2 \, R \, r$.

To begin, we'll plot the circle using the functions $x(t) = r \, \cos t$ and $y(t) = r \, \sin t$ for $t_0 \leq t \leq 2 \, \pi + t_0$, where $t_0$ is an angle such that
\begin{equation} \cos t_0 = \dfrac {A}{\sqrt{A^2 + B^2}} \qquad \text{and} \qquad \sin t_0 = \dfrac {B}{\sqrt{A^2 + B^2}}. \end{equation}

\noindent The distance from a point $\bigl( x(t), \, y(t) \bigr)$ to the line is given by the function 
\begin{equation} \begin{aligned} r(t) & = \frac {\bigl| A \, x(t) + B \, y(t) - C \bigr|}{\sqrt{A^2 + B^2}} = \bigl| \cos t_0 \cdot r \,  \cos t + \sin t_0 \cdot r \,  \sin t - R \bigr| \\ & = R - r \, \cos (t - t_0). \end{aligned} \end{equation}

\noindent (We have used the Angle Difference Formula $\cos (t - t_0) = \cos t \, \cos t_0 + \sin t \, \sin t_0$.  Recall that $R > r$ by assumption.) We have the arc length differential
\begin{equation} ds = \sqrt{ \bigl[ x'(t) \bigr]^2 + \bigl[ y'(t) \bigr]^2} \, dt  = \sqrt{ \bigl[ -r \, \sin t \bigr]^2 + \bigl[ r \, \cos t \bigr]^2} \, dt = r \, dt. \end{equation}

\noindent The desired area of revolution is given by the integral
\begin{equation} 2 \, \pi \int_{t_0}^{t_1} r(t) \, d s = 2 \, \pi \, r \int_{t_0}^{2 \pi + t_0} \bigl( R - r \, \cos (t - t_0) \bigr) \, dt = 4 \, \pi^2 \, R \, r. \end{equation}

%\section{Appendix}

%All of the graphs shown in this article were generated using \texttt{Mathematica}.  For the convenience of the reader, we include the code which generates these graphs.  The input is a line $A \, x + B \, y = C$, entered as $\texttt{\{A,B,C\}}$; two functions $x = x(t)$ and $y = y(t)$, entered as $\texttt{\{x,y\}}$; and an interval $t_0 \leq t \leq t_1$, entered as $\texttt{\{t,t0,t1\}}$.

%\bibliography{area}
%\bibliographystyle{plain}

\end{document}